\documentclass[12pt]{amsart}
\usepackage{fullpage,url,amssymb}
\usepackage[all]{xy} 

\DeclareFontEncoding{OT2}{}{} 

\usepackage{color}

\newcommand{\Pfister}[1]{\langle \langle #1 \rangle \rangle}
\newcommand{\defi}[1]{\textsf{#1}}

\newcommand{\Aff}{{\mathbb A}}
\newcommand{\C}{{\mathbb C}}
\newcommand{\F}{{\mathbb F}}

\newcommand{\N}{{\mathbb N}}

\newcommand{\Q}{{\mathbb Q}}
\newcommand{\R}{{\mathbb R}}
\newcommand{\Z}{{\mathbb Z}}
\newcommand{\Qbar}{{\overline{\Q}}}

\newcommand{\Fbar}{{\overline{\F}}}
\newcommand{\Xbar}{{\overline{X}}}

\newcommand{\pp}{{\mathfrak p}}
\newcommand{\mm}{{\mathfrak m}}

\newcommand{\calF}{{\mathcal F}}

\newcommand{\calL}{{\mathcal L}}

\newcommand{\OO}{{\mathcal O}}

\DeclareMathOperator{\Th}{Th}

\DeclareMathOperator{\trdeg}{tr deg}
\DeclareMathOperator{\Krdim}{Kr dim}

\DeclareMathOperator{\Char}{char}

\DeclareMathOperator{\End}{End}

\DeclareMathOperator{\Gal}{Gal}

\newcommand{\smooth}{{\operatorname{smooth}}}

\newcommand{\surjects}{\twoheadrightarrow}

\newcommand{\isom}{\simeq}

\newcommand{\intersect}{\cap} 
\newcommand{\tensor}{\otimes}

\newtheorem{theorem}{Theorem}[section]
\newtheorem{lemma}[theorem]{Lemma}
\newtheorem{corollary}[theorem]{Corollary}
\newtheorem{proposition}[theorem]{Proposition}

\theoremstyle{definition}
\newtheorem{definition}[theorem]{Definition}
\newtheorem{question}[theorem]{Question}

\theoremstyle{remark}
\newtheorem{remark}[theorem]{Remark}

\usepackage[alphabetic,backrefs,lite]{amsrefs} 

\begin{document}

\title[Uniform first-order definitions]{Uniform first-order definitions in finitely generated fields}
\author{Bjorn Poonen}
\thanks{This research was supported by NSF grant DMS-0301280
          and a Packard Fellowship.  The author thanks the 
        Isaac Newton Institute for hosting a visit in the summer of 2005.}
\address{Department of Mathematics, University of California, 
        Berkeley, CA 94720-3840, USA}
\email{poonen@math.berkeley.edu}
\urladdr{http://math.berkeley.edu/\~{}poonen}
\date{July 22, 2005}

\begin{abstract}
We prove that there is a first-order sentence in the language of rings
that is true for all finitely generated fields of characteristic~$0$
and false for all fields of characteristic~$>0$.
We also prove that for each $n \in \N$,
there is a first-order formula $\psi_n(x_1,\ldots,x_n)$
that when interpreted in a finitely generated field $K$
is true for elements $x_1,\ldots,x_n \in K$
if and only if the elements are algebraically dependent over 
the prime field in $K$.
\end{abstract}

\maketitle

\section{Introduction}\label{S:introduction}

\subsection{Theorems for finitely generated fields}

It is important, especially when trying to transfer results
from one structure to another, to know whether a property 
can be expressed by the truth of a first-order sentence.
For example, it is a basic theorem of model theory 
that a first-order sentence in the language of rings
true for one algebraically closed field of characteristic~$0$
is true for all algebraically closed fields of characteristic~$0$;
it is because of this that many theorems proved for $\C$ 
using analytic methods are known to hold for 
arbitrary algebraically closed fields of characteristic~$0$.

But many properties have no such simple first-order characterizations.
Characterizing characteristic~$0$ fields among all fields
with a single sentence is a typical impossible task:
compactness shows that for every sentence $\phi$ 
valid for an algebraically closed field of characteristic~$0$, 
there is a number $p_0$ such that $\phi$ holds also for every algebraically
closed field of characteristic $\ge p_0$.

Similarly, 
one cannot detect whether elements $t_1,\ldots,t_n$ 
are algebraically dependent over the prime subfield 
by means of a single formula with $n$ free variables.
When $n=1$, this would amount to defining 
the relative algebraic closure $k$ of the prime field in $K$ uniformly,
but if $K=\C$, 
then every definable subset of $K$ is either finite or cofinite, 
while $k=\Qbar$ is neither.

We focus our attention on finitely generated fields,
i.e., fields that are finitely generated over the prime subfield.
Our main theorems show that the arithmetic of these fields
is rich enough that the previously impossible tasks become possible.

Before giving these theorems, we fix a few conventions.
A \defi{formula} is a first-order formula in the language of rings.
If $\phi$ is a \defi{sentence} (a formula with no free variables)
and $K$ is a field, then $K \models \phi$ is the statement
that $\phi$ is true for $K$.
A \defi{definable subset} is one defined by a formula 
{\em without constants}.
The \defi{prime subfield} $\F$ of a field $K$ is its minimal subfield
(either $\Q$ or $\F_p$ for some prime $p$).
When discussing a finitely generated field $K$,
we will always use $k$ to denote
the \defi{field of constants},
defined as the (relative) algebraic closure of $\F$ in $K$.

\begin{theorem}
\label{T:characteristic}
There is a sentence that is true
for all finitely generated fields of characteristic~$0$,
and false for {\em all} fields of characteristic~$>0$.
\end{theorem}

\begin{theorem}
\label{T:prime field}
There exists a formula $\phi(t)$ 
that when interpreted in an {\em infinite} finitely generated field $K$
is true if and only if $t \in \F$.
\end{theorem}

The hypothesis in Theorem~\ref{T:prime field} that $K$ is infinite
is necessary, because there is no uniform definition of $\F_p$
in $\F_{p^2}$ \cite{Chatzidakis-vandenDries-Macintyre1992}.

\begin{theorem}
\label{T:constants}
There exists a formula $\psi(t)$
that when interpreted in a finitely generated field $K$
is true if and only if $t \in k$.
\end{theorem}

\begin{theorem}
\label{T:alg dep}
For each $n \in \N$,
there exists a formula $\psi_n(t_1,\ldots,t_n)$
that when interpreted in a finitely generated field $K$
is true if and only if $t_1,\ldots,t_n$ are algebraically dependent
over $k$ (or equivalently, over $\F$).
\end{theorem}

\begin{remark}
The $n=1$ case of Theorem~\ref{T:alg dep} is Theorem~\ref{T:constants}.
\end{remark}

\begin{remark}
Using different methods,
it is possible to prove geometric analogues of these results,
for function fields over algebraically closed and other ``large'' fields.
These will appear in a future joint paper with F.~Pop.
\end{remark}

\subsection{Questions about the richness of the arithmetic of finitely generated fields}

The sentence of Theorem~\ref{T:characteristic} 
defines the class of characteristic~$0$ finitely generated fields
among all finitely generated fields.
Can any ``reasonable'' class 
of finitely generated fields can be distinguished by a single sentence?

E.~Hrushovski suggested the following definition of ``reasonable''.
Fix any natural bijection
between the set of $(r,f_1,\ldots,f_m)$ with $r \in \N$
and $f_1,\ldots,f_m \in \Z[x_1,\ldots,x_r]$
and a recursive subset $A \subseteq \N$,
such as the one sending $(r,f_1,\ldots,f_m)$
to the concatenation of the ASCII symbols
of the characters in the \TeX\ code for $(r,f_1,\ldots,f_m)$.
There is an algorithm for testing whether an ideal
$(f_1,\ldots,f_m) \subseteq \Z[x_1,\ldots,x_r]$ is prime;
i.e., the set of $a \in A$ corresponding to $(r,f_1,\ldots,f_m)$
with this property is a recursive subset $B$.
We have a surjection
\[
        \kappa \colon B \to 
        \{\text{isomorphism classes of finitely generated fields}\}
\]
sending an $a \in B$ corresponding to $(r,f_1,\ldots,f_m)$
to the fraction field of $\Z[x_1,\ldots,x_r]/(f_1,\ldots,f_m)$.
Call a set $S$ 
of isomorphism classes of finitely generated fields \defi{reasonable}
if $\{a \in B: \kappa(a) \in S\}$
is a first-order definable subset in $(\N,+,\cdot)$.
If $\phi$ is a sentence, 
the set of isomorphism classes of finitely generated fields $K$
such that $K \models \phi$ is a reasonable set.

\begin{question}
\label{Q:reasonable}
Does every reasonable set of isomorphism classes of finitely generated fields
arise from a sentence $\phi$ as above?
\end{question}

For example, the set of (isomorphism classes of)
function fields of smooth projective $\Q$-varieties 
having a rational point is a reasonable set,
so a positive answer to Question~\ref{Q:reasonable} 
would say in particular that there is
a sentence that is true for these finitely generated fields
and no others.
Also, any single isomorphism class forms a reasonable set,
so a positive answer to Question~\ref{Q:reasonable}
would imply a positive answer to the following:

\begin{question}
\label{Q:single out K}
Is it true that for every finitely generated field $K$,
there is a sentence $\phi_K$
that is true for $K$ 
and false for all finitely generated fields $L \not\isom K$?
\end{question}

The \defi{theory} $\Th(K)$ of a field $K$
is the set of sentences $\phi$ such that $K \models \phi$.
Fields $K$ and $L$ are \defi{elementarily equivalent}
(and we then write $K \equiv L$)
if $\Th(K)=\Th(L)$.
A positive answer to Question~\ref{Q:single out K}
would imply a positive answer to the following question,
which was raised by G.~Sabbagh in a special case
and formulated explicitly by F.~Pop~\cite{Pop2002}.

\begin{question}
\label{Q:ee}
Is it true that 
whenever $K$ and $L$ are finitely generated fields and $K \equiv L$,
we have $K \isom L$?
\end{question}

\subsection{Structure of this paper}

Section~\ref{S:previous} describes earlier work that will be useful to us.
Section~\ref{S:char0} shows how to define
the field of constants in a finitely generated field of characteristic $0$.
Section~\ref{S:charp} shows how to define
the family of relatively algebraically closed global function fields 
in a finitely generated field of characteristic $>0$:
this is by far the hardest part of the paper.
Section~\ref{S:proofs} combines the results obtained so far
to prove the main theorems.
Finally, Section~\ref{S:negative} shows that the formulas
promised by the main theorems cannot be purely existential.

\section{Previous work}
\label{S:previous}

\begin{definition}
The \defi{Kronecker dimension} of a finitely generated field $K$ is
\[
        \Krdim K := \begin{cases} 
                \trdeg(K/\F_p), & \text{ if $\Char K=p>0$} \\
                \trdeg(K/\Q)+1, & \text{ if $\Char K=0$.}
        \end{cases}
\]
A \defi{global field} is a finitely generated field $K$
with $\Krdim K=1$;
such a field is either a \defi{number field} (finite extension of $\Q$)
or a \defi{global function field} 
(function field of a curve over a finite field).
\end{definition}

R.~Rumely \cite{Rumely1980}, building on the work of J.~Robinson,
proved Theorems \ref{T:characteristic}, \ref{T:prime field}
and \ref{T:constants} 
(and hence also Theorem~\ref{T:alg dep})
in the case of global fields.
We record some of his results in the following theorem.

\begin{theorem}\hfill
\label{T:Rumely}
\begin{enumerate}
\item \label{I:rumely-char}
There is a sentence that is true for all number fields
and false for all global function fields.
\item \label{I:rumely-prime}
There is a formula that for any global field defines the prime subfield.
\item \label{I:rumely-Z}
There is a formula that for any number field defines the subset $\Z$.
\item There is a formula that for any number field defines the subset $\N$.
\item \label{I:rumely-constants}
There is a formula that for any global function field
defines the constant field.
\item There is a formula with a parameter $x$ 
that for all global function fields defines the ring $\F[x]$ in $K$.
(See \cite{Rumely1980}*{p.~211}.)
\item There is a formula with a parameter $x$ 
that for all global function fields defines the ring $k[x]$ in $K$.
\item There are formulas with a parameter $x$ that for 
any global field $K$ and $x \in K-k$
defines a model of $(\N,+,\cdot)$ in $K$
consisting of the powers of $x$.
\end{enumerate}
\end{theorem}

We will also need the following result proved by Pop 
using the recently proved connection between isotropy of Pfister forms
and cohomological dimension: see also Proposition~\ref{P:Voevodsky}.

\begin{theorem}[\cite{Pop2002}]
\label{T:Krdim}
For each $n \in \N$,
there is a sentence $\sigma_n$ that for a finitely generated field $K$
holds if and only if $\Krdim K=n$.
\end{theorem}

It follows from these results that Question~\ref{Q:single out K}
has a positive answer for any finitely generated $K$
with $\Krdim K \le 1$.

The paper~\cite{Pop2002} also answered 
Question~\ref{Q:ee} and its geometric analogue
in the case where one of the finitely generated fields
is a function field of general type.

\section{Defining constants in characteristic $0$}
\label{S:char0}

\begin{definition}
Suppose $P$ (the parameter space) is a definable subset of $K^N$,
and $D$ is a definable subset of $K^M \times P$.
For each $\vec{p} \in P$, we get
a subset $S_{\vec{p}}:=\{\vec{a} \in K^M : (\vec{a},\vec{p}) \in D\}$.
Any family of subsets of $K^M$ that equals $\{S_{\vec{p}} : \vec{p} \in P\}$
for some such $P$ and $D$ will be called a
\defi{definable family of subsets}.
\end{definition}

If $X$ is a variety over a field $K$, and $L$ is a field extension of $K$,
we let $X_L$ denote $X \times_K L$.

%

\begin{lemma}
\label{L:EK=Ek}
Given a finitely generated field $K$ of characteristic~$0$,
with field of constants $k$,
there exists an elliptic curve $E$ over $\Q$ such that 
$E(\Q)$ is infinite and $E(K)=E(k)$.
\end{lemma}

\begin{proof}
This is a special case of \cite{Moret-Bailly2005preprint}*{Lemma~11.1(i)}.
\end{proof}

The idea behind the following two proofs is contained in
the proof of \cite{Denef1978-polynomial}*{Lemma 3.4(iii)}
and \cite{Kim-Roush1995}*{Proposition~3}:
see also \cite{Moret-Bailly2005preprint}*{Lemma~11.1(ii)}.

\begin{lemma}
\label{L:E(Q_p)}
Let $E$ be an elliptic curve over $\Q_p$.
Equip $E(\Q_p)$ with the $p$-adic topology.
The closure of any infinite subgroup of $E(\Q_p)$
is an open neighborhood of the identity $O \in E(\Q_p)$.
\end{lemma}

\begin{proof}
By the theory of formal groups,
$E(\Q_p)$ contains a finite-index subgroup
that is isomorphic as a topological group to $\Z_p$.
The result for infinite subgroups of $E(\Q_p)$ 
now follows from the corresponding result for $\Z_p$.
\end{proof}

\begin{lemma}
\label{L:family-dense}
There exists a definable family $\calF_1$ of subsets of $K$,
defined by a formula independent of $K$,
such that if $K$ is a finitely generated field of characteristic~$0$,
some $S \in \calF_1$ is a subset of $k$ such that for each finite prime $p$,
the intersection $S \cap \Q$ is $p$-adically dense in $\Q_p$.
\end{lemma}

\begin{proof}
For each $(a,b) \in K^2$, consider the subset
\[
        S_{a,b}:=\{x/y : x \in K, \; y \in K^*, \text{ and } y^2 = x^3 + a x + b \}.
\]
of $K$.
These subsets form a definable family $\calF_0$.

Now suppose that $K$ is a finitely generated field of characteristic~$0$.
Let $a,b$ be the elements of $\Q$ defining the elliptic curve $E$
of Lemma~\ref{L:EK=Ek}.
Then $E(K)=E(k)$, so $S_{a,b} \subseteq k$.
Applying Lemma~\ref{L:E(Q_p)} to the infinite group $E(\Q)$
shows that $E(\Q)$ is $p$-adically dense 
in a neighborhood of $O$ in $E(\Q_p)$.
Since the rational function $x/y$ on $E$ is a uniformizing parameter at $O$,
it is a local diffeomorphism between a neighborhood of $O$ in $E(\Q_p)$
and a neighborhood of $0$ in $\Q_p$.
Thus the $p$-adic closure of $S_{a,b} \cap \Q$ contains a neighborhood 
of $0$ in $\Q_p$.

Finally, for each $(a,b)$, let $T_{a,b}$ be the set of ratios
of nonzero elements of $S_{a,b}$.
The subsets $T_{a,b}$ form a definable family $\calF_1$
having the required property.
\end{proof}

\begin{lemma}
\label{L:Pfister}
Let $K$ be a finitely generated field of characteristic~$0$
with $\sqrt{-1} \in K$.
Let $S$ be as in Lemma~\ref{L:family-dense}.
Then $k$ is the set of $t \in K$ such that
for all $s_1,s_2,s_3 \in S$,
the Pfister form $\Pfister{s_1,s_2,t-s_3}$ over $K$ represents $0$.
\end{lemma}

\begin{proof}
Suppose $t \in k$ and $s_1,s_2,s_3 \in S$.
Then $s_1,s_2,t-s_3 \in k$.
The form $\Pfister{s_1,s_2,t-s_3}$ is of rank $>4$,
and $k$ is a totally complex number field,
so $\Pfister{s_1,s_2,t-s_3}$ represents $0$ over $k$,
and hence also over $K$.

Now suppose $t \notin k$.
Let $V$ be an integral $k$-variety with function field $K$.
Replacing $V$ by an open subset, 
we may assume $t$ is regular on $V$.
Since $t \notin k$, the map $t\colon V \to \Aff^1_k$ is dominant.
Since $\Char k=0$, we may shrink $V$ further
to assume that $t$ is smooth.
Since $t$ is dominant, some $s_3 \in S \subseteq \Aff^1(k)$ 
is in the image of $t$.
Choose a closed point $v \in V$ at which $t-s_3$ vanishes.
Then $t-s_3$ is part of a system of local parameters at $v$.
Let $\ell$ be the residue field of $v$, so $\ell$ is a number field.
Choose a prime $p$ of $\Q$ that splits completely in $\ell$.
Choose $a,b \in \Q_p$ such that $\Pfister{a,b}$ has no
nontrivial zero.
Approximate $a,b$ $p$-adically by $s_1,s_2 \in S \cap \Q$,
closely enough that $\Pfister{s_1,s_2}$ 
still has no nontrivial zero over $\Q_p$.
Since $\ell$ injects into $\Q_p$,
$\Pfister{s_1,s_2}$ has no nontrivial zero over $\ell$.
Lemma~\ref{L:local parameters} implies that $\Pfister{s_1,s_2,t-s_3}$ 
has no nontrivial zero over $K$.
\end{proof}

\begin{lemma}
\label{L:k}
There exists a definable family $\calF_2$ of subsets of~$K$,
defined by a formula independent of~$K$,
such that if $K$ is a finitely generated field of characteristic~$0$,
one of the subsets in $\calF_2$ is the field of constants~$k$.
\end{lemma}

\begin{proof}
Let $\calF_1$ be as in Lemma~\ref{L:family-dense}.
For each $S \in \calF_1$, consider 
\[
        A_S:=\{t \in K: (\forall s_1,s_2,s_3 \in S)\; \Pfister{s_1,s_2,t-s_3} \text{ represents $0$ over $K(\sqrt{-1})$}\}.
\]
The family $\calF_2$ of such subsets is a definable family
defined by a formula independent of $K$.

Now suppose $K$ is a finitely generated field of characteristic~$0$.
If $S$ is the infinite subset of $k$ promised by Lemma~\ref{L:family-dense},
then Lemma~\ref{L:Pfister} implies that 
\[
        A_S = K \intersect \text{(constant subfield of $K(\sqrt{-1})$)} = k.
\]
\end{proof}

\begin{lemma}
\label{L:defining constants}
There is a formula $\phi(x)$
that when interpreted in a finitely generated field
of characteristic~$0$ defines its field of constants,
and that when interpreted in a finitely generated field
of characteristic~$>0$ defines the empty set.
\end{lemma}

\begin{proof}
Let $\calF_2$ be as in Lemma~\ref{L:k}, so $k \in \calF_2$.
We can find a definable subfamily $\calF_3 \subseteq \calF_2$
consisting of the sets in $\calF_2$ that are fields.
All these fields are finitely generated,
so we may apply Theorem~\ref{T:Krdim} with $n=1$
and then Theorem~\ref{T:Rumely}(\ref{I:rumely-char})
to find a definable subfamily $\calF_4 \subseteq \calF_3$
consisting of the sets of $\calF_3$ that are number fields.
The union $\ell$ of the sets in $\calF_4$ is a subset of $K$
definable by a formula independent of $K$.

If $K$ is a finitely generated field of characteristic~$0$,
each set in $\calF_4$ is a subfield of $k$,
and $k \in \calF_4$,
so $\ell=k$.
On the other hand, if $K$ is a finitely generated field of 
characteristic~$p>0$, then $\calF_4$ is empty, so $\ell=\emptyset$.
\end{proof}

\section{Defining subfields of transcendence degree $1$ over finite fields}
\label{S:charp}

For this section, we fix the following notation:
\begin{itemize}
\item $K$ is a finitely generated field of characteristic $p>0$.
\item $k$ is the field of constants of $K$.
\item $\calL$ is the set of subfields $L \subseteq K$ 
such that $L$ is relatively algebraically closed in $K$
and $\trdeg(L/k)=1$.
\end{itemize}
Moreover, we use the following notation for the rest of the paper:
let $d=3$ if $\Char k=2$ and $d=2$ otherwise.

Our main goal in this section is to construct a definable family of subsets
of $K$ consisting of the fields in $\calL$.
We will mainly try to follow the constructions of the previous section
used to define the collection of number fields in $K$ 
in the characteristic~$0$ case, 
but more elaborate arguments will be needed, since certain steps
fail.
For instance, the analogue of Lemma~\ref{L:E(Q_p)}
for local fields of characteristic~$p$ is false.

\subsection{Nearly prime generalized Mersenne numbers}

Let $E$ be an ordinary elliptic curve over a finite field $\F_q$.
We would like to find primes $\ell$ such that $\#E(\F_{q^\ell})$
is close to being prime.
Let $R=\End_{\F_q} E$.
Let $F \in R$ be the $q$-power Frobenius endomorphism.
Let $\alpha,\beta \in \Qbar$ 
be the eigenvalues of $F$ acting on a Tate module.
Then $\#E(\F_{q^\ell})=(\alpha^\ell-1)(\beta^\ell-1)$
(whereas Mersenne numbers are numbers of the form $2^\ell-1$).

To quantify this, for an integer $n$, define $\Psi(n):=\sum_{p|n} \frac{1}{p}$,
where the sum is over the distinct prime divisors $p$ of $n$.
We consider $n$ to be nearly prime if $\Psi(n)$ is small.

\begin{lemma}
\label{L:Mersenne}
For any ordinary elliptic curve $E$ over $\F_q$,
\[
        \liminf_{\ell \to \infty} 
        \Psi\left( \frac{\#E(\F_{q^\ell})}{\#E(\F_q)} \right) = 0.
\]
\end{lemma}

\begin{proof}
Let $\OO$ be the ring of integers of $R \tensor \Q = \Q(\alpha)$.
Define 
\[
        e_\ell:=\frac{\#E(\F_{q^\ell})}{\#E(\F_q)} 
        = N_{F/\Q}\left( \frac{\alpha^\ell-1}{\alpha-1} \right).
\]
If $\ell$ and $m$ are distinct primes,
then the ideal $(x^\ell-1,x^m-1)$ in $\Z[x]$ equals $(x-1)$,
since reducing the smaller generator modulo the other iteratively
amounts to running the Euclidean algorithm on the exponents.
Thus the ideals $\left( \frac{\alpha^\ell-1}{\alpha-1} \right)$
and $\left( \frac{\alpha^m-1}{\alpha-1} \right)$
of $\OO$ are coprime.
There are at most prime ideals of $\OO$ above a given prime of $\Z$,
so each prime of $\Z$ divides at most two of the $e_\ell$.
Hence
\begin{equation}
\label{E:mersenne1}
        \sum_{\ell \le B} \Psi(e_\ell) 
        \le 2 \Psi \left(\prod_{\ell \le B} e_\ell \right).
\end{equation}
Since $e_\ell \le \#E(\F_{q^\ell}) = q^{\ell}(1+o(1))$ as $\ell \to \infty$, 
we have $\prod_{\ell \le B} e_\ell = q^{O(B^2)}$ as $B \to \infty$.
In particular $\prod_{\ell \le B} e_\ell$ has at most $O(B^2)$ 
distinct prime factors.
Thus, if $p_j$ is the $j$-th prime, 
\begin{equation}
\label{E:mersenne2}
        \Psi\left(\prod_{\ell \le B} e_\ell \right) \le \sum_{j=1}^{O(B^2)} \frac{1}{p_j}.
\end{equation}
The prime number theorem implies $p_j = (1+o(1)) j \log j$ as $j \to \infty$,
so 
\begin{equation}
\label{E:mersenne3}
        \sum_{j=1}^{O(B^2)} \frac{1}{p_j} = O\left(\sum_{j=2}^{O(B^2)} \frac{1}{j \log j} \right) = O\left( \int_{2}^{O(B^2)} \frac{dx}{x \log x} \right) = O(\log \log O(B^2)) = O(\log \log B).
\end{equation}
Combining equations 
\eqref{E:mersenne1}, \eqref{E:mersenne2}, \eqref{E:mersenne3}
gives 
\[
        \sum_{\ell \le B} \Psi(e_\ell) = O(\log \log B).
\]
But there are $\frac{B}{\log B}(1+o(1))$ primes $\ell$ up to $B$,
so some term $\Psi(e_\ell)$ on the left is bounded by
\[
        O\left(\frac{\log \log B}{B/\log B}\right) = 
        O\left(\frac{(\log B )(\log \log B)}{B}\right),
\]
which tends to $0$ as $B \to \infty$.
\end{proof}

\begin{remark}
The conclusion of Lemma~\ref{L:Mersenne}
holds also for supersingular $E$,
but we do not need this.
\end{remark}

\begin{lemma}
\label{L:probability}
Let $E$ be an ordinary elliptic curve over $\F_q$.
Fix $\ell \in \Z_{\ge 1}$.
The probability that a random element of $E(\F_{q^\ell})$
generates $E(\F_{q^\ell})/E(\F_q)$ as an $R$-module
is at least $1-2\Psi(n)$ where $n:=\#E(\F_{q^\ell})/\#E(\F_q)$.
\end{lemma}

\begin{proof}
Since $E$ is ordinary, \cite{Lenstra1996}*{Theorem~1(a)}
gives an isomorphism of $R$-modules
$\iota\colon E(\F_{q^\ell}) \to R/(F^\ell-1)R$.
If $P \in E(\F_{q^\ell})$ does not generate $E(\F_{q^\ell})/E(\F_q)$ 
as an $R$-module,
then the $R$-submodule of $E(\F_{q^\ell})$ generated by $P$ and $E(\F_q)$
must correspond under $\iota$ to a proper ideal of $R':=R/(F^\ell-1)R$,
and hence be contained in a maximal ideal $\mm$ of $R'$
of residue characteristic dividing $n$.
If $\mm_1,\ldots,\mm_r$ are the maximal ideals of residue characteristic~$p$,
then $\mm_1 \intersect \cdots \intersect \mm_r$ 
is of index at least $p^r$ in $R'$,
but it contains $pR'$ which has index at most $(R:pR)=p^2$,
so $r \le 2$.
The probability that $\iota(P)$ lies in a given $\mm$
of residue characteristic~$p$ is at most $1/p$,
so the probability that $\iota(P)$ lies in some $\mm$
of residue characteristic dividing $n$
is at most $\sum_{p|n} \frac{2}{p} = 2 \Psi(n)$.
\end{proof}

\subsection{More elliptic curve lemmas}

\begin{lemma}
\label{L:ordinary}
For any finite field $k$,
there exists an ordinary elliptic curve over $k$.
\end{lemma}

\begin{proof}
Honda-Tate theory~\cite{Honda1968} gives $E$ over $k$ 
for which the trace of Frobenius equals $1$.
\end{proof}

\begin{lemma}
\label{L:sum of z-coordinates}
There is a universal constant $c \in \R_{>0}$ such that the following holds.
Let $E$ be an elliptic curve over a finite field $\F_q$
in Weierstrass form.
Let $z$ be the rational function $y/x$ on $E$.
Let $G$ be a subgroup of $E(\F_q)$.
If $(E(\F_q):G) < c q^{1/2}$,
\[
        \{z(g_1)+z(g_2): g_i \in G - \{\text{poles of $z$}\}\} = \F_q.
\]
\end{lemma}

\begin{proof}
Define a curve
\[
        X:= \{(P_1,P_2) \in (E-\{\text{poles of $z$}\})^2 : 
                                z(P_1)+z(P_2) = t\}.
\]
Let $\Xbar$ be the Zariski closure of $X$ in $E \times E$.
Let $X^\smooth$ be the smooth locus of $X$,
and $\Xbar^\smooth$ the smooth locus of $\Xbar$.
If we use $w=1/z$ as a uniformizer at the identity $O$ of $E$,
then we obtain a system of local parameters $w_1,w_2$
at $(O,O) \in E \times E$,
and the local equation of $X$ there is $1/w_1+1/w_2=t$,
which after clearing denominators is $w_2+w_1=t w_1 w_2$;
thus $(O,O) \in \Xbar^\smooth$.
Let $g \ge 1$ be a universal upper bound 
(independent of $q$, $E$, $t$) for the geometric genus
of the geometric components of the normalization $\hat{X}$ of $\Xbar$,
and let $r$ be a universal upper bound for the number of 
$\Fbar_q$-points in $\hat{X} - X^\smooth$.
Define $c := 1/(2g+r)$.

The separable isogeny $F-1\colon E \to E$ 
factors through $E \surjects E':=E/G$
so we get the homomorphism $\phi$ in the commutative diagram
\[
\xymatrix{
E \ar[rr]^{F-1} \ar@{->>}[d] && E \ar@{->>}[d] \\
E' \ar@{.>}[urr]^\phi \ar[rr]_{F-1} && E'. \\
}
\]
The lower right triangle now shows that $G=\phi(E'(\F_q))$.
Let $\delta:=\deg \phi = (E(\F_q):G) \le c q^{1/2}$.

Let $X'=(\phi,\phi)^{-1}(X) \subset E' \times E'$,
and define $\Xbar'$ and $(\Xbar')^\smooth$ similarly as inverse images.
Then $(O,O) \in (\Xbar')^\smooth(\F_q)$,
so the irreducible component of $\Xbar'$ containing $(O,O)$
is geometrically irreducible.
Thus the dense open subset ${X'}^\smooth$ of $\Xbar'$ also
has a geometrically irreducible component $C'$.
Since $C'$ is finite \'etale of degree at most $\delta$ over $X^\smooth$,
it is a curve of geometric genus at most 
$g':=\delta(g-1)+1 \le c q^{1/2} g$ 
with at most $r':= \delta r \le c q^{1/2} r$ geometric points removed.
By the Weil conjectures 
\[
        \#C'(\F_q) \ge q + 1 -2 g' q^{1/2} - r' > q - 2 c g q - c r q^{1/2}  \ge q(1-(2g+r)c) = 0.
\]
In particular, there exists a point $(P_1,P_2) \in C'(\F_q)$.
For $i=1,2$, define $g_i:=\phi(P_i) \in \phi(E'(\F_q)) = G$.
By definition of $C'$, we have $z(g_1)+z(g_2)=t$.
\end{proof}

\begin{lemma}
\label{L:prime degree points}
Let $E$ be an ordinary elliptic curve over a finite field $\F_q$.
Let $R=\End_{\F_q}(E)$.
Let $W$ be a geometrically integral $\F_q$-variety
with a smooth morphism $\pi=(\pi_1,\pi_2)\colon W \to \Aff^1 \times E$.
For some prime $\ell>3$,
there exists $w \in W(\F_{q^\ell})$
such that $\pi_1(w) = z(Q_1) + z(Q_2)$
for some $Q_1,Q_2 \in R \cdot \pi_2(w)$.
\end{lemma}

\begin{proof}
Let $S_1:=\pi(W(\F_{q^\ell}))$.
Let $S_2$ be the set of $(t,P) \in (\Aff^1 \times E)(\F_{q^\ell})$
such that $t=z(Q_1) + z(Q_2)$
for some $Q_1,Q_2 \in R \cdot P$.
We need to show that $S_1$ and $S_2$ intersect,
so it will suffice to prove that
\begin{equation}
\label{E:needed}
        \#S_1 + \#S_2 > \#(\Aff^1 \times E)(\F_{q^\ell})
\end{equation}
for some prime $\ell>3$.

By the Weil conjectures,
$\#W(\F_{q^\ell}) = q^{\ell(\dim W)}(1+o(1))$ as $\ell \to \infty$.
Let $b$ be a bound on the number of components 
of the geometric fibers of $W \to \Aff^1 \times E$;
then each fiber of $W(\F_{q^\ell}) \to (\Aff^1 \times E)(\F_{q^\ell})$
has size at most $(b+o(1)) q^{\ell(\dim W - 2)}$.
Dividing, we obtain
\begin{equation}
\label{E:S1}
        \#S_1 \ge \left(\frac{1}{b} - o(1) \right) q^{2 \ell}.
\end{equation}

By Lemma~\ref{L:Mersenne},
we can find infinitely many primes $\ell$ such that
the integer $n:=\#E(\F_{q^\ell})/\#E(\F_q)$ satisfies 
$\Psi(n) \le \frac{1}{4b}$;
we assume from now on that $\ell$ satisfies this.
By Lemma~\ref{L:probability},
the probability that a random $P\in E(\F_{q^\ell})$
generates $E(\F_{q^\ell})/E(\F_q)$ as an $R$-module
is at least $1-2\Psi(n) \ge 1-\frac{1}{2b}$.
In this case, $G:=R\cdot P$ has index at most $\#E(\F_q)$
in $E(\F_{q^\ell})$, 
so for $\ell \gg 1$, 
Lemma~\ref{L:sum of z-coordinates} applied to $E$ over $\F_{q^\ell}$ 
implies that
\[
        \{z(g_1)+z(g_2): g_i \in G - \{\text{poles of $z$}\}\} = \F_{q^\ell}.
\]
Thus
\begin{equation}
\label{E:S2}
        \#S_2 \ge \left(1-\frac{1}{2b} \right) \#E(\F_{q^\ell}) \cdot \# \F_{q^\ell} = \left(1 -\frac{1}{2b} - o(1) \right) q^{2\ell}.
\end{equation}

Finally,
\begin{equation}
\label{E:S3}
        \#(\Aff^1 \times E)(\F_{q^\ell}) = \left(1 + o(1) \right) q^{2\ell}.
\end{equation}
Combining \eqref{E:S1}, \eqref{E:S2}, and~\eqref{E:S3}
gives \eqref{E:needed} if $\ell$ is large enough.
\end{proof}

Suppose $L \in \calL$.
Let $E$ be an elliptic curve over $k$ in Weierstrass form.
For $u \in L-k$, define a degree-$2$ \'etale $L$-algebra
$L_u:= L \tensor_{k(u)} k(E)$
where the homomorphism $k(u) \to k(E)$ sends $u$ to the 
coordinate function $x \in k(E)$ 
for a fixed Weierstrass model,
and let $E_u$ be the elliptic curve over $L$
obtained by twisting $E_L$ by $L_u/L$.
Define $K_u:= K \tensor_{k(u)} k(E)$ similarly.
(If $\Char k \ne 2$ and the Weierstrass model is $y^2=f(x)$,
then $E_u$ is $f(u) y^2 = f(x)$.)


\begin{lemma}
\label{L:good twist}
For any $L \in \calL$ and elliptic curve $E$ over $k$,
there exists $u \in L-k$ such that $L_u$ and $K_u$ are fields
and $E_u(K)=E_u(L)$.
\end{lemma}

\begin{proof}
Let $V$ be an integral $L$-variety with function field $K$.
Let $A$ be the Albanese variety of $V$.
Let $A_1,\ldots,A_n$ be the distinct $L$-simple abelian varieties
appearing in a decomposition of $A$ up to isogeny.
Choose a nontrivial place $v$ of $L$ at which all the $A_i$ 
have good reduction, and choose $u$ so that $v(u)=-1$.
Since $v(u)$ is negative and odd,
$v$ totally ramifies in $L_u/L$;
in particular $L_u$ is a field.
Since $L$ is relatively algebraically closed in $K$,
the algebra $K_u$ is a field too.

Now $E_u$ is a twist of $E$ by a quadratic extension $L_u/L$ 
in which $v$ ramifies, 
so $E_u$ has bad reduction at $v$.
Hence $E_u$ is not isogenous to any $A_i$.
Therefore every morphism $A \to E_u$ is constant.
So every $L$-rational map $V \dashrightarrow E_u$ is constant.
In other words, $E_u(K)=E_u(L)$.
\end{proof}

Suppose $u$ is as in Lemma~\ref{L:good twist}.
Let $\sigma$ be the nontrivial element of $\Gal(K_u/K)$.
We may identify $E_u(K)$ with 
\[
        E(K_u)^{\sigma=-1}:=\{P \in E(K_u): {}^\sigma P = -P\}.
\]

\subsection{Defining subfields of transcendence degree~$1$}

\begin{lemma}
\label{L:anisotropy}
Suppose $L \in \calL$.
Suppose that $E$ is an ordinary elliptic curve over $k$.
Choose $u \in L-k$ as in Lemma~\ref{L:good twist}.
Let $U_1=x(E_u(K)) - \{0\}$.
Let $U_2$ be the set of $u_2 \in K$
expressible as $z(Q_1)+z(Q_2)$
for some $Q_1,Q_2 \in E(K_u) - \{\text{poles of $z$}\}$
with $x(Q_1),x(Q_2) \in U_1$.
Suppose $c \in k^{\times} - k^{\times d}$.
For $t \in K$,
we have $t \in L$ if and only if
for all $u_1 \in U_1$ and $u_2 \in U_2$
\[
        \Pfister{t-u_2,u_1^{-1}}_d \tensor \langle 1,-c \rangle_d
\]
has a nontrivial zero over $K$.
\end{lemma}

\begin{proof}

Suppose $t \in L$.
By choice of $u$, we have $E_u(K)=E_u(L)$,
so $U_1 \subseteq L$.
The conditions $x(Q_1),x(Q_2) \in U_1$
imply $Q_1,Q_2 \in L_u$,
so each $u_2 \in U_2$ is algebraic over $L$,
and hence in $L$.
Then Proposition~\ref{P:Voevodsky}
implies that $\Pfister{t-u_2,u_1^{-1}}_d \tensor \langle 1,-c \rangle_d$
has a nontrivial zero over $L$,
and hence also over $K$.

Now suppose $t \in K-L$.
Thus $t$ is transcendental over $L$.
View $k(E)$ as a subfield of $L_u$.
Then $k(E)(t)$ is the function field of $\Aff^1 \times E$ over $k$.
Let $M_0 \subseteq K_u$ be a purely transcendental extension of $k(E)(t)$ 
such that $[K_u:M_0]<\infty$.
Let $M$ be the relative separable closure of $M_0$ in $K_u$.
Since $L_u/L$ is ramified somewhere,
each of $k$, $L_u$, $K_u$ is relatively algebraically closed in the next,
so $k$ is relatively algebraically closed in $M$.
So there is a geometrically integral $k$-variety $W$
with function field $M$.
Shrinking $W$,
we may assume that the rational map $W \dashrightarrow \Aff^1 \times E$
corresponding to the extension $M$ over $k(E)(t)$
is a smooth morphism $\pi=(\pi_1,\pi_2)$
whose image is disjoint from $\Aff^1 \times E[2]$.
Let $\ell$, $w$, $Q_1$, $Q_2$
be as provided by Lemma~\ref{L:prime degree points} (with $\F_q=k$).
Let $(\tau,P)=\pi(w) \in (\Aff^1 \times E)(\F_{q^\ell})$.
Thus $\tau = z(Q_1)+z(Q_2)$,
and for $i=1,2$
we have $Q_i=\rho_i(P)$ for some $\rho_i \in R:=\End_k(E)$.
Let $\eta$ be a separable endomorphism of $E$ killing $P$:
for instance, we could take $\eta=F^\ell-1$.
The point $\bar{\eta} \in E(k(E))$ corresponding to $\eta$
is in $E(K_u)^{\sigma=-1}$,
so $u_1:=x(\bar{\eta}) \in k(E) \subseteq K_u$
belongs to $U_1$.
Let $\bar{\rho}_i$ be the point of $E(k(E))$ corresponding to $\rho_i$.
Let $u_2=z(\bar{\rho}_1)+z(\bar{\rho}_2) \in U_2$.
Since $P \notin E[2]$,
the rational function $u_1$ on $E$ has a simple pole at $P$.
The value of the rational function $t-u_2$ on $\Aff^1 \times E$
at $(\tau,P)$ is $\tau - (z(\rho_1(P)) + z(\rho_2(P))) = 0$.
Thus $t-u_2,u_1^{-1}$ 
are a system of local parameters at $(\tau,P)$ on $\Aff^1 \times E$.
We pull them back to $W$.
Since $W \to \Aff^1 \times E$ is smooth,
we can extend this pair to a system of local parameters at $w$ on $W$.
Since the residue field of $w$ is $\F_{q^\ell}$ with $\ell$ prime to $d$,
the form $\langle 1,-c \rangle_d$ has no nontrivial zero over $\F_{q^\ell}$.
By Lemma~\ref{L:local parameters}, the form
\[
        Q:=\Pfister{t-u_2,u_1^{-1}}_d \tensor \langle 1,-c \rangle_d
\]
has no nontrivial zero over $M$.
By Corollary~\ref{C:purely inseparable},
$Q$ has no nontrivial zero over $K_u$,
so it has no nontrivial zero over $K$.
\end{proof}

\begin{lemma}
\label{L:family containing L}
There exists a definable family of subsets of $K$,
defined by a formula independent of $K$,
such that if $K$ is a finitely generated field of characteristic $p>0$
and $K$ contains a primitive cube root of $1$ if $p=2$,
then all $L \in \calL$ belong to the family.
\end{lemma}

\begin{proof}
For each elliptic curve $E$ over $K$ and $u,c \in K$,
let $U_1$ and $U_2$ be as in Lemma~\ref{L:anisotropy},
and define $S_{E,u,c}$ to be the set of $t \in K$
such that for all $u_1 \in U_1$ and $u_2 \in U_2$,
\[
        \Pfister{t-u_2,u_1^{-1}}_d \tensor \langle 1,-c \rangle_d
\]
has a nontrivial zero over $K$.
We can quantify over $E$
by quantifying over the vector of coefficients 
$\vec{a}=(a_1,a_2,a_3,a_4,a_6)$ of a Weierstrass equation
subject to the nonsingularity constraint $\Delta(a_1,\ldots,a_6) \ne 0$.
Thus $\{S_{E,u,c}\}$ is a definable family.

Given $L \in \calL$,
choose an ordinary $E$ over $k$ as in Lemma~\ref{L:ordinary},
and choose $u \in L-k$ as in Lemma~\ref{L:good twist}.
Also we can choose $c \in k^{\times} - k^{\times d}$
since either $\#k$ is odd and $d=2$,
or $\#k$ is even and $d=3$ and $k$ contains a primitive cube root of~$1$.
Then $S_{E,u,c}=L$ by Lemma~\ref{L:anisotropy}.
\end{proof}

\begin{proposition}
\label{P:defining all L}
There exists a definable family of subsets of $K$,
defined by a formula independent of $K$,
such that if $K$ is a finitely generated field of characteristic $p>0$,
the subsets in the family are exactly those in $\calL$.
\end{proposition}

\begin{proof}
Let $K'=K(\zeta)$ where $\zeta^2+\zeta+1=0$.
Take the family of Lemma~\ref{L:family containing L}
for $K'$, and intersect each set with $K$.
Using Theorem~\ref{T:Krdim},
retain only the sets in the family that are fields of Kronecker dimension~$1$.
\end{proof}

\section{Proofs of theorems}
\label{S:proofs}

\begin{proof}[Proof of Theorem~\ref{T:characteristic}]
Using Lemma~\ref{L:defining constants} 
and Theorem~\ref{T:Rumely}(\ref{I:rumely-Z}),
we can find a formula defining a set $S$
such that whenever $K$ is finitely generated and of characteristic~$0$,
we get $S=\Z$.
Use the sentence that says that $S$ is closed under addition
and $S \ne 2S$ and $2 \ne 0$.
This sentence is true for $S=\Z$,
but false for any subset of a field of positive characteristic.
\end{proof}

Now that Theorem~\ref{T:characteristic} is proved,
we may handle the characteristic~$0$ and $>0$ cases separately in proving 
Theorems \ref{T:prime field}, \ref{T:constants}, and~\ref{T:alg dep}.

\begin{proof}[Proof of Theorem~\ref{T:prime field}]
If $\Char K=0$, combine Lemma~\ref{L:defining constants}
with Theorem~\ref{T:Rumely}(\ref{I:rumely-prime}).
If $\Char K>0$, use Theorem~\ref{T:Rumely}(\ref{I:rumely-prime})
to take the prime subfield of each member of the 
family given by Proposition~\ref{P:defining all L},
and take their intersection:
this works as long as the family is nonempty,
which holds since $K$ is infinite.
\end{proof}

\begin{proof}[Proof of Theorem~\ref{T:constants}]
If $\Char K=0$, use Lemma~\ref{L:defining constants}.
If $\Char K>0$, the field $k$
is the set of elements belonging
to the field of constants of each relatively algebraically closed
subfield of transcendence degree~$1$:
this is uniformly definable
by Theorem~\ref{T:Rumely}(\ref{I:rumely-constants})
and Proposition~\ref{P:defining all L}.
\end{proof}

\begin{remark}
Using Proposition~\ref{P:defining all L}
and the theorems just proved,
we can also extend other results in~\cite{Rumely1980}
to the finitely generated case.
For instance, 
there is a formula $\phi(x,y)$ such that
when $K$ is an infinite finitely generated field of characteristic~$>0$
and $x \in K$,
we have $\{y \in K: \phi(x,y)\} = \F[x]$.
The same can be done for $k[x]$ in place of $\F[x]$.
\end{remark}

Most of the remaining work in proving Theorem~\ref{T:alg dep}
is contained in the following lemma,
whose proof is very close to that of Fact~1.3(3) in~\cite{Pop2002}.

\begin{lemma}
\label{L:alg dep over global fields}
For each $n \in \N$, there exists a formula $\phi_n(t_1,\ldots,t_n)$
in the language of fields augmented by a predicate for a subfield
such that the following holds.
Let $K$ be a finitely generated extension of a global field $L$,
and assume $L$ is relatively algebraically closed in $K$.
Then elements $t_1,\ldots,t_n \in K$ are algebraically dependent over $L$
if and only if $\phi_n(t_1,\ldots,t_n)$ holds over $K$
with the predicate corresponding to $L$.
\end{lemma}

\begin{proof}
The elements $t_1,\ldots,t_n$ are algebraically dependent over $L$
if and only if they are algebraically dependent over $L(\sqrt{-1})$,
so by replacing $K,L$ by $K(\sqrt{-1}),L(\sqrt{-1})$,
we may reduce to the case that $\sqrt{-1} \in L$.
Similarly we may assume that $L$ contains $\zeta$
satisfying $\zeta^2+\zeta+1=0$.

It will suffice to show that
$t_1,\ldots,t_n$ are algebraically dependent over $L$
if and only if
for all $a,b,c_1,\ldots,c_n \in L$,
the form
\[
        q:=\Pfister{t_1-c_1,\ldots,t_n-c_n,a}_d \tensor \langle 1,-b \rangle_d
\]
has a nontrivial zero.

If $t_1,\ldots,t_n$ are algebraically dependent over $L$
and $a,b,c_1,\ldots,c_n \in L$,
then $q$ has a nontrivial zero by Proposition~\ref{P:Voevodsky}.

Now suppose $t_1,\ldots,t_n$ are algebraically independent over $L$.
Extend to a transcendence basis $t_1,\ldots,t_N$ of $K$ over $L$.
Let $M$ be the maximal separable extension of $L(t_1,\ldots,t_N)$ in $K$.
Choose a smooth integral $L$-variety $W$ with function field $M$.
Shrinking $W$, we may assume that 
$\tau:=(t_1,\ldots,t_N)\colon W \to \Aff^N_L$
is an \'etale morphism.
Choose $c:=(c_1,\ldots,c_N) \in \Aff^N(L)$ in the image of $\tau$.
and let $w \in W$ be a closed point in the fiber $\tau^{-1}(c)$.
Then $t_1-c_1,\ldots,t_N-c_N$ are a system of local parameters at $w$ on $W$.
Let $L'$ be the residue field of $w$, so $L'$ is a finite extension of $L$.
Choose a nonarchimedean place $\pp$ of $L$ that splits completely in $L'$,
so $L'$ injects into the completion $L_\pp$.
Choose a $\pp$-adic unit $b \in L$ whose residue is not a $d$-th power,
and choose $a \in L$ of $\pp$-adic valuation $1$.
By Lemma~\ref{L:DVR}, 
$\Pfister{a}_d \tensor \langle 1,-b \rangle_d$ 
has no nontrivial zero over $L_\pp \supset L'$.
By Lemma~\ref{L:local parameters},
$q$ has no nontrivial zero over $M$.
By Corollary~\ref{C:purely inseparable},
it also has no nontrivial zero over $K$.
\end{proof}

\begin{proof}[Proof of Theorem~\ref{T:alg dep}]
We may assume $n \ge 1$.
If $\Char K=0$, Lemma~\ref{L:alg dep over global fields} does the job.
If $\Char K>0$, 
we have that $t_1,\ldots,t_n$ are algebraically dependent
if either $t_1 \in k$
or there exists $L \in \calL$ such that $t_1 \in L$ 
and $t_2,\ldots,t_n$ are algebraically dependent over $L$:
this is uniformly definable by a formula,
by Theorem~\ref{T:constants},
Proposition~\ref{P:defining all L},
and Lemma~\ref{L:alg dep over global fields}.
\end{proof}

\section{Negative results for existential definitions}
\label{S:negative}

In this section, 
we show that existential formulas cannot satisfy the requirements
of Theorems \ref{T:characteristic},
\ref{T:prime field}, 
\ref{T:constants}, 
\ref{T:alg dep}.

Given an existential formula,
we can convert each polynomial inequality $f(x_1,\ldots,x_n) \ne 0$
to $(\exists y) f(x_1,\ldots,x_n)y=1$
and convert each disjunction of polynomial equalities $f=0 \lor g=0$
to $fg=0$.
Thus we need only consider formulas given as a conjunction of
polynomial equalities, preceded by existential quantifiers.

The following gives a negative result for Theorem~\ref{T:characteristic}.

\begin{proposition}
There is no existential sentence that is true for $\Q$
and false for $\F_p$ for all primes $p$.
\end{proposition}

\begin{proof}
If a closed subscheme $V$ of $\Aff^n_\Z$ has a $\Q$-point $P$,
then it has an $\F_p$-point for any $p$ not dividing the denominators
of the coordinates of $P$.
\end{proof}

The following gives a negative result for Theorem~\ref{T:prime field}.

\begin{proposition}
\label{P:negative prime field}
There is no existential formula that defines the prime field $\F$
in every infinite finitely generated field $K$.
\end{proposition}

\begin{proof}
Suppose such a formula exists.
Then there is a closed subscheme $V$ of $\Aff^n_\Z$
such that if $\pi\colon V \to \Aff^1_\Z$ is the projection onto the
first coordinate,
$\pi(V(K))=\F$
for every infinite finitely generated field $K$.
The morphism $\pi$ must be dominant,
since otherwise $\#\pi(V(K))$ would be bounded for all finitely generated $K$
of sufficiently large characteristic, while $\#\F$ is unbounded.

Choose an irreducible component $V_0$ of $V$ 
that dominates $\Aff^1_\Z$.
Take $K$ to be the function field of $V_0$.
The value of $\pi$ at the tautological point of $V_0(K)$
is not in $\F$, contradicting the assumption on $V$.
\end{proof}

The following gives a negative result for Theorem~\ref{T:constants},
and hence also for the more general Theorem~\ref{T:alg dep}.

\begin{proposition}
For each fixed $p \ge 0$,
there is no existential formula that defines the field of constants
for all finitely generated fields of characteristic $p$.
\end{proposition}

\begin{proof}
Repeat the proof of Proposition~\ref{P:negative prime field},
observing that the size of the field of constants is unbounded,
even when we fix the characteristic.
\end{proof}

\appendix
\section{Diagonal forms}

Let $d$ be a positive integer.
and let $a_1,\ldots,a_n$ be elements of a field $K$.
Then $\langle a_1,\ldots,a_n \rangle_d$ denotes
the diagonal form
\[
        a_1 x_1^d + \cdots + a_n x_n^d \in K[x_1,\ldots,x_n].
\]
Define the tensor product of two such forms 
$\langle a_1,\ldots,a_m \rangle_d$ 
and
$\langle b_1,\ldots,b_n \rangle_d$ 
to be the diagonal form in $mn$ variables 
whose coefficients are the products $a_i b_j$.
Finally, define
\begin{align*}
        \Pfister{a}_d 
        &:= \langle 1,a,\ldots,a^{d-1} \rangle_d \\
        \Pfister{a_1,\ldots,a_n}_d 
        &:= \Pfister{a_1}_d \tensor \cdots \tensor \Pfister{a_n}_d,
\end{align*}
so $\Pfister{a_1,\ldots,a_n}_d$ is a diagonal degree-$d$ form
in $d^n$ variables.
If $d=2$, then $\Pfister{a_1,\ldots,a_n}_d$ 
is called a \defi{Pfister form}.

\begin{proposition}
\label{P:Springer}
Let $q(x_1,\ldots,x_n)$ be a homogeneous form over a field $K$,
and let $L$ be a finite extension of $K$.
Suppose that either $\deg q=2$ and $[L:K]$ is odd,
or $\deg q=3$ and $[L:K]=2$.
If $q$ has a nontrivial zero over $L$,
then $q$ has a nontrivial zero over $K$.
\end{proposition}

\begin{proof}
This is well known: see \cite{LangAlgebra}*{Chapter~V, Exercise~28}.
\end{proof}

\begin{corollary}
\label{C:purely inseparable}
Let $K$ be a field.
Let $d=3$ if $\Char K = 2$ and $d=2$ otherwise.
Let $q(x_1,\ldots,x_n)$ be a homogeneous form of degree $d$ over a field $K$.
Let $L$ be a purely inseparable extension of $K$.
If $q$ has a nontrivial zero over $L$,
then $q$ has a nontrivial zero over $K$.
\end{corollary}

\begin{proof}
If $q$ has a nontrivial zero over $L$,
the coordinates of this zero generate a finite purely inseparable extension
of $K$.
By induction, we reduce to the case $[L:K]=p$, where $p:=\Char K$.
Now the result follows from Proposition~\ref{P:Springer}.
\end{proof}

\begin{proposition}
\label{P:Voevodsky}
Let $k$ be a separably closed field,
a finite field, or a number field;
define $\epsilon$ to be $0$, $1$, or $2$, respectively.
Let $d=3$ if $\Char k = 2$ and $d=2$ otherwise.
If $k$ is a number field, assume that $\sqrt{-1} \in k$.
Let $K$ be a finitely generated extension of $k$
of transcendence degree $r$.
If $n \ge r+\epsilon$ and $m \ge 2$ 
then for all $a_1,\ldots,a_n,b_1,\ldots,b_m \in K$,
the form
\[
        \Pfister{a_1,\ldots,a_n}_d \tensor \langle b_1,\ldots,b_m \rangle_d
\]
has a nontrivial zero over $K$.
\end{proposition}

\begin{proof}
The separably closed case reduces to the algebraically closed case
by Corollary~\ref{C:purely inseparable}.
If $k$ is algebraically closed or finite,
then $k$ is a $C_\epsilon$ field in the sense of \cite{Lang1952},
and $K$ is a $C_{r+\epsilon}$ field,
so the result follows.
If $k$ is a number field, use \cite{Pop2002}*{Fact~1.3(1)}.
\end{proof}

\begin{lemma}
\label{L:DVR}
Let $K$ be a field with discrete valuation $v\colon K^\times \surjects \Z$.
Let $\OO$ be the valuation ring, let $\pi \in K$ be such that $v(\pi)=1$,
and let $k=\OO/(\pi)$.
Let $d \in \Z_{\ge 2}$.
Let $q$ be a diagonal degree-$d$ form over $\OO$
whose reduction modulo $\pi$ has no nontrivial zero over $k$.
Then the form $q':= q \tensor \Pfister{\pi}_d$ 
has no nontrivial zero over $K$.
\end{lemma}

\begin{proof}
Write
\[
        q' = q(\vec{x}_0) + \pi q(\vec{x}_1) + \cdots 
                + \pi^{d-1} q(\vec{x}_{d-1}).
\]
If the coordinates of $\vec{x}_0$ are in $\OO$
and not all in $\pi\OO$,
then $v(q(\vec{x}_0))=0$, since $q$ has no nontrivial zero in $k$.
More generally, if $\vec{x}_0$ is nonzero, it is a power of $t$ times such
a primitive vector, so $v(q(\vec{x}_0)) \equiv 0 \pmod{d}$.
Similarly, if $\vec{x}_i$ is nonzero,
then $v(\pi^i q(\vec{x}_i)) \equiv i \pmod{d}$.
Since these valuations are distinct (when not $+\infty$),
the form $q'$ has no nontrivial zero over $K$.
\end{proof}

The following is close to results used in \cite{Pop2002}.

\begin{lemma}
\label{L:local parameters}
Let $k$ be a field,
and let $V$ be an integral $k$-variety with function field $K$.
Suppose that $v$ is a smooth closed point on $V$,
and that $t_1,\ldots,t_m$ are part of a system of local parameters at $v$.
Let $d \in \Z_{\ge 2}$.
Let $q$ be a diagonal degree-$d$ form over $k$ 
having no nontrivial zero the residue field of $v$.
Then $q \tensor \Pfister{t_1,\ldots,t_m}_d$ 
has no nontrivial zero over $K$.
\end{lemma}

\begin{proof}
We may assume that $t_1,\ldots,t_m$ is a complete 
system of local parameters (i.e., $m=\dim X$).
Let $k'$ be the residue field of $v$.
Since $v$ is a smooth point,
$K$ embeds into the iterated Laurent series field $k'((t_1))\cdots((t_m))$.

We prove by induction that for $j=0,\ldots,m$,
the form $q \tensor \Pfister{t_1,\ldots,t_j}_d$
has no nontrivial zero over $k'((t_1))\cdots((t_j))$.
The case $j=0$ is given,
and Lemma~\ref{L:DVR} provides the inductive step.

Taking $j=m$ gives the result.
\end{proof}

\section*{Acknowledgements} 

I thank Ehud Hrushovski for suggesting a precise definition
of ``reasonable'' class of finitely generated fields.

\end{document}